\documentclass[a4paper,12pt]{article}
\usepackage{amssymb,amsmath,amscd,amsfonts,amsthm}
\usepackage{latexsym,graphpap}
\usepackage[dvips]{epsfig}
\usepackage{color}
\usepackage[english,activeacute]{babel}

\newtheorem{theorem}{Theorem}
\newtheorem{prop}[theorem]{Proposition}

\newenvironment{demo}{ \noindent \emph{\textbf{Proof:}}}{\hfill$\square$\\}

\newcommand{\RR}{\mathbb{R}}

\newcommand{\Cc}{\mathcal{C}}
\newcommand{\Nc}{\mathcal{N}}
\newcommand{\Ac}{\mathcal{A}}
\newcommand{\Bc}{\mathcal{B}}
\newcommand{\Uc}{\mathcal{U}}

\newcommand{\grad}{\nabla}
\newcommand{\Un}{1\hspace{-1.5mm}1}
\newcommand{\no}{n$^{\text{o}}$}

\newcommand\nor[2]{\left\|#1\right\|_{#2}}
\newcommand{\pc}{ \usefont{T1}{cmtl}{m}{n} \selectfont}

\newdimen\texpscorrection
\newdimen\figcenter
\def\figurewithtex #1 #2 #3 #4 #5\cr{\null
  {\goodbreak\figcenter=\hsize\relax
  \advance\figcenter by -#4truecm
  \divide\figcenter by 2
  \begin{figure}[hbt]
  \vskip #3truecm\noindent\hskip\figcenter
  \includegraphics{#1}{\hskip\texpscorrection\input #2 }
  \vskip 0.8truecm{\baselineskip=0.8\baselineskip
  \noindent \vbox{\noindent {\footnotesize #5}}\par}
  \end{figure}}}
\def\point#1 #2 #3 {\rlap{\kern #1 truecm
\raise #2 truecm \hbox{#3}}}

\numberwithin{equation}{section}

\begin{document}

\title{\bf A note on the {semi-}global controllability of the semilinear
wave equation}

\author{ Romain \textsc{Joly}\footnote{Institut Fourier - UMR5582
CNRS/Universit\'e de Grenoble - 100, rue des Maths - BP 74 - F-38402
St-Martin-d'H\`eres, France, email: {\pc romain.joly@ujf-grenoble.fr}} {~\&~}
Camille \textsc{Laurent}\footnote{CNRS, UMR 7598, Laboratoire Jacques-Louis
Lions, F-75005, Paris, France  }
\footnote{UPMC Univ Paris 06, UMR 7598, Laboratoire Jacques-Louis Lions,
F-75005, Paris, France, email: {\pc laurent@ann.jussieu.fr }}}

\maketitle
\vspace{1cm}

\begin{abstract}
We study the internal controllability of the semilinear wave equation  
$$v_{tt}(x,t)-\Delta v(x,t) + f(x,v(x,t))= \Un_{\omega} u(x,t)$$
for some nonlinearities $f$ which can produce several non-trivial steady states.

One of the usual hypotheses to get {semi-}global controllability, is to
assume that
$f(x,v)v\geq 0$. In this case, a stabilisation term $u=\gamma(x)v_t$ makes any
solution converging to zero. The {semi-}global controllability then follows
from a
theorem of local controllability and the time reversibility of the 
equation. 

In this paper, the nonlinearity $f$ can be more general, so that 
the solutions of the damped equation may converge to another equilibrium than
$0$. To prove {semi-}global controllability, we study the controllability
inside a
compact attractor and show that it is possible to travel from one equilibrium
point to another by using the heteroclinic orbits.\\[2mm]
\end{abstract}

\section{Introduction}\label{sect-one}

{Let $\Omega$ be a smooth connected riemannian manifold of dimension $d$
with boundary and let $\omega$ be an open subset of $\Omega$. Let
$X=H^1_0(\Omega)\times L^2(\Omega)$, let $\Delta$ be the Laplacian operator
with Dirichlet boundary conditions and let
$f\in\Cc^{1}(\overline\Omega\times\RR,\RR)$.
We consider the internal control of the wave equation
\begin{equation}\label{eq}
\left\{\begin{array} {ll} 
v_{tt}(x,t)-\Delta v(x,t) + f(x,v(x,t))= \Un_{\omega}
u(x,t)~~~~~&(x,t)\in\Omega\times (0,T)\\
v(x,t)=0 & (x,t)\in\partial\Omega\times (0,T)\\
(v,\partial_t v)(x,0)=V_0 \in X &       \end{array}\right.
\end{equation}}
where the function $u\in L^1((0,T),L^2(\Omega))$ is the control and
$V=(v,v_t)\in\Cc^0([0,T],X)$ is the state of the system. 

This paper concerns the {semi-}global controllability of the
semilinear wave equation, that is that we are interested in proving the
following property:
\begin{itemize}
 \item[{\bf (SGC)}] For any bounded subset $\Bc$ of
$X$, there exists a time $T(\Bc)>0$ such that, for any
$V_0$ and $V_1$ in $\Bc$, there exist $T\leq T(\Bc)$ and $u\in
L^1((0,T),L^2(\Omega))$ such that the solution of \eqref{eq} satisfies
$V(T)=(v,v_t)(T)=V_1$.
\end{itemize}

Other articles consider the boundary control problem by
replacing \eqref{eq} with a control from the boundary for the same
semilinear equation. We will not study these boundary conditions while our
methods could apply, but it would require some further analysis of the
non-homogeneous boundary value problem for the nonlinear equation. 

{
The control theory for the linear case $f\equiv 0$ is now well
known. The almost necessary and sufficient condition for global controllability
is that $\omega$ satisfies the geometric control condition: there exists $L>0$
such that any generalized geodesic of
$\Omega$ of length $L$ meets the set $\omega$ where the control is effective.
Note that the condition is always sufficient, see the works of
Bardos, Lebeau and Rauch \cite{BLRinterne} and \cite{BLR}. However, it is
necessary only if we replace $1_{\omega}(x)$ by a smooth function
$\gamma_{\omega}(x)$ satisfying $\omega= \left\{x\left|\gamma_{\omega}(x)\neq
0\right.\right\}$, see Burq-G\'erard \cite{BurqGerardCNS} for a proof of the
necessity for the closely related boundary control. In the case of a control
with  $1_{\omega}$, some subtleties can happen if one ray stays close to
$\partial\omega$, as for example, for the control from a half hemisphere of
$S^2$, see \cite{LebeauControlHyp} p.174.}

The study of the semilinear case $f\neq 0$ is still in progress (see
\cite{Coron} and \cite{ZZ}).
To our knowledge, the control problem is solved mainly in three different cases:
\begin{itemize}
\item {\it Local control.} We assume {$f(x,0)=0$}, that is that $v\equiv 0$
is a steady state of \eqref{eq}.
The problem of local controllability states as property (SGC), except that $\Bc$
is not any bounded set
but is a small neighborhood of $0$. The local controllability is known for both
internal and boundary control, see \cite{Chewning} and \cite{Zuazua-loc}. 

The result of Coron and Tr\'elat \cite{Coron-Trelat} is quite different but
could be described roughly as local controllability near a whole connected
component of the set of steady states. The authors study local controllability
by boundary control near a path of steady states, in dimension $d=1$. Notice
that their method using quasi-static deformations, enables to compute
effectively the control. Remark also that an important part of our proof
consists in proving that we can go from an equilibrium to another using the
global compact attractor and is therefore similar, in spirit, to their result.

\item {\it Quasilinear nonlinearity.} Another class of problems for which (SGC)
holds, is the case where $f$ is almost linear, for instance globally
lipschitzian {or super-linear but with a growth of the type $s\ln^\beta
s$},
see {\cite{FYZ},} \cite{LT}, \cite{Li-Zhang}, \cite{Zuazua-control}, 
\cite{Coron} {and the references therein. Notice that these papers prove
in fact the global controllability, that is that the time $T$ in (SGC) can be
chosen independently of $\Bc$.}

\item {\it Sign condition {$f(x,s)s\geq 0$}.} In \cite{Zuazua-stab},
\cite{DLZ} and \cite{RomCam}, stabilization results are proved, that is that, if
{$f(x,s)s\geq 0$} and if $f$ satisfies some additional conditions, all the
solutions of the damped wave equation
\begin{equation}\label{eq-d-intro}
 v_{tt}(x,t)-\Delta v(x,t) + {f(x,v(x,t))}= -\Un_{\omega} v_t(x,t)
\end{equation}
goes to $0$ when $t$ goes to $+\infty$. Adding a local control result near $0$,
one easily gets a global control result by using $u=-v_t$ as control in
\eqref{eq}. Therefore, \cite{Zuazua-stab}, \cite{DLZ} and \cite{RomCam}
also yield the {semi-}global controllability of semilinear wave equations,
under a                   
sign assumption for $f$. Notice that this sign condition cannot be removed
carelessly if $f$ is super-linear, due to the counter-example of
\cite{Zuazua-control}.
\end{itemize}

The purpose of this paper is to release the sign condition {$f(x,s)s\geq
0$} to an asymptotic sign condition and to {allow more complicated
dynamics for the damped wave equation \eqref{eq-d-intro}}. In this case, the
solutions of \eqref{eq-d-intro} may converge to another equilibrium than $0$.
Therefore, one has to explain how to move from one equilibrium to another in the
control problem \eqref{eq}.

{
At this point, we would like to make a remark about the dependence of $f$ on
the space variable $x$. In fact, in the references given above, $f$ is assumed
to be independent of $x$, in order to simplify the calculations and because it
seems to be in the habit to do so in control theory. On the opposite, from
the dynamical point of view, richer dynamics are
welcome and it is in the habit to allow $f$ to depend on $x$. We will allow 
this $x-$dependence in this paper, however, we underline that:
\begin{list}{--}{\itemsep=0pt \leftmargin=5mm}
\item the results of the references cited above should also hold for $f\equiv
f(x,v)$,
as soon as the $x-$dependence does not change the important properties of $f$,
as growth estimates, sign conditions\ldots
\item the result and the discussions of this paper are also meaningful for
$f\equiv f(v)$, since functions as $f(v)=\lambda v(v-1)(v+1)$ can also generate
non-trivial dynamics for \eqref{eq-d-intro}. 
\end{list}
}

\section{Main result}

{The main idea of this paper is very general and may be applied to several
kind
of problems. However, to avoid too abstract formalisms, we state our main result
in the following particular cases, corresponding to some already known results,
that we will use.

We assume that one of these sets of assumptions is satisfied:}
\begin{itemize}
\item[\bf Case A.] $\Omega$ is a smooth bounded open domain of $\RR^d$ and
there exists $x_0\in\RR^d$ such that 
\begin{equation}\label{hyp-omega}
\{ x\in\partial\Omega ~/~(x-x_0).\nu >0\}\subset \overline\omega~.
\end{equation}
Moreover, if
$d\geq 2$, we assume that there exists $0\leq p <d/(d-2)$ and $C>0$ such that 
\begin{equation}\label{hyp-f2}
|f'_x(x,s)|+|f(x,s)|\leq C(1+|s|)^p~~\text{ and }~~|f'_s(x,s)|\leq
C(1+|s|)^{p-1}~.
\end{equation}
\begin{equation}\label{hyp-f1}
\inf_{x\in\overline\Omega} ~ \liminf_{|s|\rightarrow \infty} \frac{f(x,s)}s > 
\lambda_1~,
\end{equation} where $\lambda_1$ is the first (non-positive) eigenvalue of
$\Delta$.
\item[\bf Case B.] $\Omega=\RR^3$ and $\omega$ is the exterior of a ball.
Moreover, $f$ satisfies \eqref{hyp-f2} with $p \in [0,5)$, and for 
any $(x,s)\in\Omega\times\RR$, 
\begin{equation}\label{hyp-f3}
\left(x\not\in B(0,R) \text{ or }|s|\geq R\right)~~\Longrightarrow~~ f(x,s)s\geq
cs^2~.
\end{equation}
for a positive constant $c>0$.
\item[\bf Case C.] $\Omega$ is a smooth compact manifold with boundary of
dimension $d=3$. We assume that the geodesics of $\overline{\Omega}$ do not have
contact of infinite order with $\partial \Omega$ and that 
$\omega$ satisfies the geometric control condition of \cite{BLR}. Moreover,
$f(x,s)$ is {of class $\Cc^\infty$,} analytic with respect to $s$, and
satisfies \eqref{hyp-f2} with $p\in[0,5)$. We also ask
\begin{equation}\label{hyp-f4}
|s|\geq R~~\Longrightarrow~~ f(x,s)s > 0~.
\end{equation}
\end{itemize}

\begin{theorem}\label{th}
Assume that the above hypotheses, in one of the cases A, B or C, hold. Then,
the wave equation \eqref{eq} is {semi-}globally controllable in the sense
that, for any
bounded subset $\Bc$ of $X$, there exists $T(\Bc)>0$ such that, for any $V_0$
and $V_1$ in $\Bc$, there exist $T\leq T(\Bc)$ and $u\in L^1((0,T),L^2(\Omega))$
such that the solution of \eqref{eq} satisfies $V(T)=(v,v_t)(T)=V_1$.
\end{theorem}

As said in Section \ref{sect-one}, Theorem \ref{th} was already known when
\eqref{hyp-f1} is replaced by the more restricted sign condition
$f(x,s)s>0$ for all $s\in\RR^*$. In this simpler case, all the solutions
of the associated damped wave equation 
\begin{equation}\label{eq-damped2}
v_{tt}+\Un_\omega v_t =\Delta v - f(x,v)
\end{equation}
converge to zero. When releasing the sign assumption to the asymptotic one
\eqref{hyp-f1}, \eqref{hyp-f3} or \eqref{hyp-f4}, the dynamics of
\eqref{eq-damped2} become more complicated.
The main new idea of this paper is to show how one can use
the heteroclinic connections of \eqref{eq-damped2} to travel from one
equilibrium to another and still obtain global controllability.

The three cases A, B and C may seem restrictive, but were chosen to fit to
the already known results on existence of global compact attractors and local
control:
\begin{itemize}
\item Case A is the classical case where the local control can be proved by
multiplier methods and the Cauchy problem solved by Sobolev embedding.
\item Case B is the euclidean case with a subcritical nonlinearity which
requires Strichartz estimates. The corresponding stabilisation problem has
been studied for instance by Dehman, Lebeau and Zuazua in \cite{DLZ}.
\item Case C assumes the optimal Geometric Control Condition.
Using the global strategy of \cite{DLZ}, we have recently proved the
stabilisation and the existence of compact global attractor in this case (see
\cite{RomCam}). The Unique Continuation Property was quite tricky and required
some dynamical system tools associated with a unique continuation theorem
requiring partial analyticity. That is why we also require $f$ to be analytic
with respect to $s$ in this case. Notice that we do not have to assume that
$f(x,0)=0$ on the boundary as explained in the revised version of \cite{RomCam}
available online (see \cite{Hale-Raugel} p. 1111 for the original idea).
\end{itemize}
Of course, similar results should be true in different geometric situations.
Also, the results in case B and C could certainly be applied in any dimension.
The subcritical exponents are then $p\in [0,(d+2)/(d-2))$ for $d\geq 3$ and any
finite $p\geq 0$ for $d=1,2$.   
In fact, the method for proving Theorem \ref{th} may apply in a larger
framework, including for example boundary control. This framework is written
more explicitly in the next section.

\section{The qualitative framework of Theorem \ref{th}}

As already said, the framework behind the proof of Theorem \ref{th} is more
general than Cases A, B and C. In this section, we describe this framework
through qualitative assumptions and we show that it holds in particular for
Cases A, B and C.

\subsection{Cauchy problem} 
The assumptions of Theorem \ref{th} are such
that the wave equation \eqref{eq} satisfies the following
properties. The Cauchy problem is locally well posed and if the solution $V$ of
\eqref{eq} exists and is bounded for all $t\in [0,T]$, then there is a
neighborhood of $V_0$ such that all the solutions starting in this neighborhood
are well defined for all $t\in[0,T]$ and the Cauchy problem is continuous with
respect to the initial data. Notice that \eqref{eq} is reversible in
time and, therefore, that all these properties also hold for the Cauchy problem
backward in time. 

More precisely, in cases A, B and C, we define $X_T$ the functional spaces where
the equation will be well posed:
\begin{itemize}
\item Case A: $X_T=C([0,T],H^1_0(\Omega))\cap C^1([0,T],L^2(\Omega))$.
\item Case B and C: $X_T=C([0,T],H^1_0(\Omega))\cap C^1([0,T],L^2(\Omega))\cap 
L^{\frac{2p}{p-3}}([0,T],L^{2p}(\Omega)$
\end{itemize}
In case B and C, we need an additional Strichartz space which is chosen here to
ensure $f(x,v)\in L^1([0,T],L^2)$ when $f$ satisfies \eqref{hyp-f2} with
$p$ in the specified range (we can also assume without loss of generality that $p>3$).  
\begin{theorem}[Cauchy problem]\label{th-Cauchy1}$~$\\
For any $u\in L^1([0,T], L^2)$ and $V_0 \in
X=H^1_0(\Omega)\times L^2(\Omega)$, there exists a unique solution $v\in
X_T$ of the controlled equation \eqref{eq}. Moreover, the flow map
\begin{eqnarray*}\begin{array}{rcll}
\Phi :&X\times L^1([0,T],L^2)&\longmapsto& X_T\\
&(V_0,u)&\longrightarrow& v
\end{array} \end{eqnarray*}
is locally Lipschitz.
\end{theorem}

This local Cauchy theory is well known and the dependence of $f$ on $x$ does not
change the proof. In Case A, $f$ generates a lipschitz map in the bounded sets
of $X$ and thus all these properties of the Cauchy problem follow from the
classical semigroup theory (see \cite{Pazy}). In Cases B and C, $f$ has a higher
growth rate and Strichartz estimates are required for the local existence. 
For case B, Strichartz estimates are true in their full range and can be found
for instance in Ginibre-Velo \cite{Gi-Ve-3}. For case C, Strichartz
estimates are more recent and we refer to Burq-Lebeau-Planchon \cite{BLP} for
the first result and Blair-Smith-Sogge \cite{StrichartzSoggewave} for some wider
range of exponents. 
We refer to the books of Sogge \cite{SoggebookNLW}, Tao \cite{Tao} or Cazenave
\cite{Cazenave} (for the closely related nonlinear Schr\"odinger equation) for 
references on the Cauchy problems using Strichartz estimates. See also the
articles of Ginibre and Velo \cite{Gi-Ve-1,Gi-Ve-2}. 

{
The globalisation of the solution is obtained as follows. The natural energy of
the nonlinear wave
equation is defined by
\begin{equation}\label{energy}
E(V)=\int_\Omega \frac 12 |v_t|^2+\frac 12 |\grad v|^2 + F(x,v)~dx
\end{equation}
where $F(x,v)=\int_0^v f(x,\zeta)d\zeta$. Notice that, in any case A, B or C,
the energy $E$ is well defined because of \eqref{hyp-f2}. Bounding the
variation of the energy for a solution $V(t)$ of \eqref{eq} gives}
\begin{equation}\label{eq33}
 E(V(s))\leq E(V(t))+C\nor{u}{L^1([s,t],L^2)}\sup_{\tau \in[s,t]}E(V(\tau)).
\end{equation}
This shows that the energy remains bounded on bounded intervals. We
then use the different assumptions on the sign of $f$ to infer similar estimate
for the $H^1\times L^2$ norm. Indeed, \eqref{hyp-f1}, \eqref{hyp-f3} or
\eqref{hyp-f4} implies that $F(x,v)-\lambda_1|v|^2$ is bounded from below
everywhere and even non-negative for $x$ outside a ball $B(0,R)$. Thus, 
using Poincar\'e inequality $ |\lambda_1|\nor{v}{L^2}^2\leq \nor{\nabla
v}{L^2}^2$, we obtain the existence of $\eta>0$ such that
$$\forall V\in X~,~~E(V)\geq \eta\|V\|_X^2 + vol(B(x_0,R)) \inf  
F~.$$
Thus, the control of the energy given by \eqref{eq33} implies a bound on the
norm of the solution of \eqref{eq-damped}. Moreover, the Sobolev embeddings
$H^1(\Omega)\hookrightarrow L^{p+1}(\Omega)$ shows that the bounded sets of $X$
have a bounded energy.

\subsection{Local controllability near equilibrium points}

We will use as a black-box the local controllability near equilibrium points.
It is precisely stated as follows.
\begin{prop} \label{proplocalcontrol}
For any equilibrium $e\in H^1_0(\Omega)$ of \eqref{eq}, there exists a
neighborhood
$\Nc(e)$ of $(e,0)$ in $X$ such that \eqref{eq} is controllable in $\Nc(e)$. In 
other word, there exists a time $T(e)$ such that for any $V_0$ and $V_1$ in
$\Nc(e)$, there exist $T\leq T(e)$ and $u\in L^\infty((0,T),L^2(\Omega))$ such
that the solution of \eqref{eq} satisfies $V(T)=V_1$.
\end{prop}
 Actually, in cases A, B or C, $T(e)$ can be chosen uniformly equal to the time
of the geometric control condition. The proof of the local
controllability consists in writing the problem as a perturbation of the linear
controllability 
$$h_{tt}-\Delta h + f'_v(x,e(x))h= \Un_{\omega} u(x,t) $$
and to apply a fixed point theorem. The proof is mutatis mutandis given
in \cite{Zuazua-loc}. See also \cite{Coron}, Theorem 3 of \cite{DLZ}, Theorem
3.2 of \cite{LaurentNLW} for control close to $0$. See also
\cite{LaurentNLSdim3} where local control near trajectories are constructed for
the nonlinear Schr\"odinger equation. The only difference is that we are
close to a non trivial solution and we apply similar local control method to
$r(t,x)=v(t,x)-e(x)$ solution of 
\begin{equation}
\left\{\begin{array} {ll} 
r_{tt}(x,t)-\Delta r(x,t) + f'_v(x,e(x))r(x,t)+g(x,r)= \Un_{\omega}
u(x,t)~~~~~&(x,t)\in\Omega\times (0,T)\\
r(x,t)=0 & (x,t)\in\partial\Omega\times (0,T)\\
(r,\partial_t r)(x,0)=R_0 \in X &       \end{array}\right.
\end{equation}
where $g(x,r)=r^2\int_0^1 (1-s)f''_v(x,r(t,x)s+e(x))~ds$ is of order $2$ in $r$.

Note that for the controllability of the linear system, the geometric control
condition satisfied by $\omega$ and the fact that $f'_v(x,e(x))$ is smooth and
independent on the time are crucial for the unique continuation, see Theorem
3.8 and Corollary 4.10 of \cite{BLR}. The smoothness of $e$, and thus the one
of the potential $f'_v(x,e(x))$, follows from elliptic estimates and the
subcriticality of $f$.

\subsection{The damped wave equation is a gradient dissipative dynamical
system} 
We set $\gamma=\Un_\omega$ and we consider the damped wave equation associated
to \eqref{eq}
\begin{equation}\label{eq-damped}
\left\{\begin{array} {ll} 
v_{tt}+\gamma(x)v_t =\Delta v - f(x,v)~~~~~~~&(x,t)\in\Omega\times \RR_+~.\\
(v,v_t)(0)=V_0\in X &
\end{array}\right.
\end{equation}
{We are interested in the dynamical properties of Equation
\eqref{eq-damped}, which are recalled below. We do not give detailed proofs,
 since these results are classical or straightforward adaptations of already
existing results.  
\begin{theorem}[Cauchy problem]\label{th-Cauchy2}$~$\\
Let the assumptions of Case A, B or C be fulfilled. Then, for any $V_0\in
X=H^1_0(\Omega)\times L^2(\Omega)$ there exists a unique solution $v\in X_T$ of
the subcritical damped wave equation \eqref{eq-damped}. Moreover, this solution
is defined for all $t\in \RR$.
\end{theorem}
\begin{demo}
The local Cauchy theory follows from Theorem \ref{th-Cauchy1}. Estimating the
variation of the energy $E$, defined by \eqref{energy}, gives
\begin{equation*}
 E(V(s))\leq E(V(t))+C\int_s^t E(V(\tau))
d\tau.
\end{equation*}
When combined with Gronwall inequality, it gives backward and forward global
well-posedness, as already discussed after Theorem \ref{th-Cauchy1}. 
\end{demo}

Since \eqref{eq-damped} is globally well posed, it generates a global dynamical
system $S(t)$ on $X$, defined by $S(t)V_0=V(t)$. An important dynamical
feature is that $S(t)$ is a gradient dynamical system.
\begin{theorem}[Gradient property]\label{th-gradient}$~$\\
The energy $E$ defined by \eqref{energy} is a strict Lyapounov functional for
the dynamical system $S(t)$ generated by the damped wave equation
\eqref{eq-damped}, that is that:
\begin{enumerate}
 \item the energy $E(V(t))$ is non-increasing in time for all solutions $V(t)$
of \eqref{eq-damped},
\item if $E(V(t))$ is constant for any $t\geq 0$, then $V(t)$ is
an equilibrium point of \eqref{eq-damped}.
\end{enumerate}
\end{theorem}
\begin{demo} 
The non-increase of the energy comes from the direct computation $\frac
d{dt}E(V(t))=-\int_\omega |v_t|^2$. To show the second property, the
classical argument consists in noticing that, if $E(V(t))$ is constant, $w=v_t$
satisfies
\begin{equation}\label{eq-grad}
w_{tt}=\Delta w - f'_v(x,v)w~~~~\text{ and }~~~~~w_{|\omega}\equiv 0~.
\end{equation}
Then one uses a unique continuation result to show that $w$ vanishes
identically. In Cases A and B, using \eqref{hyp-omega}, this last property
is a straightforward consequence of a unique continuation property of \cite{KK}.
As shown in \cite{RomCam}, Case C is more involved since the geometric control
condition of \cite{BLR} is more general than (\ref{hyp-omega}). The goal is to
use the unique continuation result of \cite{Rob-Zui} to \eqref{eq-grad} to show
that $E$ is a strict Lyapounov function. This unique continuation result is
optimal in the geometric point of view, but it requires that $t\mapsto
f'(\cdot,v(\cdot,t))$ is analytic, which can be shown for analytic $f$ via an
asymptotic regularization result of \cite{Hale-Raugel}.
\end{demo}

The main result of this section is the following. We refer to \cite{Hale} and
\cite{Raugel} for introductions to the notion of attractor and gradient
dynamical system.
\begin{theorem}[Existence of a compact global attractor]\label{th-dyn}$~$\\
The dynamical system $S(t)$ generated by the damped wave
equation \eqref{eq-damped} admits a compact global attractor $\Ac$,
that is a connected compact invariant
subset of $X$, which attracts all the bounded sets of $X$. The attractor $\Ac$
consists in all the trajectories of \eqref{eq-damped}, which are globally
defined and bounded for all $t\in\RR$. 

Moreover, when $t$ goes to
$+\infty$, any trajectory of $S(t)$ converges  to the set of equilibrium points.
If the trajectory is also defined and uniformly bounded for all $t\leq 0$, then
the convergence toward the set of equilibrium points also holds for $t$ going to
$-\infty$. As a consequence, the compact global
attractor consists exactly in the set of equilibrium points and trajectories
connecting two parts of this set.
\end{theorem}
\begin{demo}
The existence of a compact global
attractor for $S(t)$ follows from the classical arguments of Theorem 3.8.5 of
Hale's book \cite{Hale} (see also Theorem 4.6 of \cite{Raugel}). 
The complete proofs can be found in \cite{Raugel} (Theorem
4.38 and the discussion above) for Case A and in Theorem 1.4 of \cite{RomCam} 
for Case C. The proof is not explicitly written in case B, but one can simply
follow the arguments given in \cite{RomCam}, using the results of \cite{DLZ}.
The three important properties are: $S(t)$ is gradient, it is asymptotically
smooth and the set of equilibrium points is bounded. The asymptotic smoothness 
comes from compactness properties of the nonlinearity $f$: in Case A,
$V\mapsto (0,f(x,v))$ is compact in $X$ due to \eqref{hyp-f2} and in Case B and
C, Hypothesis \eqref{hyp-f2} also implies some compactness property for the
nonlinearity $f$ (see Theorem 8 of \cite{BLR} and Proposition 4.3 of
\cite{RomCam}). 

Once the existence of the global attractor $\Ac$ is proved, its properties
stated in Theorem \ref{th-dyn} are the classical properties of the attractor of
a gradient dynamical system. The
convergence of the trajectories of the attractor to the set of equilibrium
points, when $t$ goes to $+\infty$ and when $t$ goes to $-\infty$, is a
consequence of Lasalle's principle (see Lemma 3.8.2 of \cite{Hale} or
Proposition 4.2 of \cite{Raugel}). The fact that the compact global attractor is
connected is a general fact and can be found for example in Lemma 2.4.1 of
\cite{Hale} of Proposition 2.19 of \cite{Raugel}. 
\end{demo}}

\noindent {\bf Remarks:}\\
{$\bullet$ For a good and simple insight into the dynamics of the attractor
$\Ac$, one may think of $\Ac$ as a finite number of isolated equilibrium points
connected by heteroclinic orbits, that is trajectories $u(t)$ converging, when
$t$ goes to $+\infty$ and $-\infty$, to two different equilibria $e_+$ and
$e_-$. Indeed, this structure is generic and typical of gradient dynamical
systems (see Theorem 3.8.5 of \cite{Hale}). However, one should be aware that
more complicated structure may occur in general, as a continuum of equilibria
and trajectories limiting to this continuum.\\}
$\bullet$ The classical case where the sign condition
$f(x,s)s\geq
0$ is assumed for all
$s\in\RR$, as in \cite{Zuazua-stab,DLZ,RomCam} for example, corresponds to the
case where the
compact global attractor $\Ac$ is reduced to $\{0\}$. Hence, in this
paper, we show how to deal with asymptotic dynamics, which are more
involved than a single equilibrium point. The main idea is to use the
heteroclinic connections to travel from an equilibrium point to another in the
attractor $\Ac$.


\section{Proof of Theorem \ref{th}}

As already mentioned, the main ideas of this paper are more general than
the framework of Theorem \ref{th}. To underline this fact, we prove our theorem
by using only the qualitative properties stated in the previous section. Hence,
our method easily extends to any control problem for which these properties
hold. 

In what follows, we will often write: "{\it we go from $V_0$ to
$V_1$}" or "{\it there exists a control bringing $V_0$ to $V_1$}" or similar
formulations. In more precise way, this will mean ``{\it there exists $T\geq 0$
and $u\in L^1([0,T],L^2)$ such that the unique solution $v$ of (\ref{eq})
satisfies $(v,\partial_t v)(T)=V_1$}''. 
In particular, the time is not fixed and can be eventually large. The dependence
of this time $T$ on the data will be evaluated after.
We call reachable state from $V_0$, and denote $\mathcal{R}(V_0)$, the set of
$V_1$ such that we can go from $V_0$ to $V_1$ (without precision on the time).

The proof of our main result can be split into several steps. The first ones
are basic facts. Step 4 consists in the key ``double U-turn'' argument. Steps
5-7 contain typical arguments coming from the study of dynamical systems. To
lighten the arguments of the proof, we will only consider the question of the
uniformity of the time of control in the end this section.

\medskip

\noindent {\bf Basic fact 1:} by Proposition \ref{proplocalcontrol}, \eqref{eq} is locally controllable
in a neighborhood $\Nc(e)$ of an equilibrium point $(e,0)$. In the following,
we assume without loss of generality that $\Nc(e)$ is a ball of $X$.

\medskip

\noindent {\bf Basic fact 2:} if $V(t)=(v,v_t)(t)$ is a solution of the damped
wave equation \eqref{eq-damped}, then, by using the control
$u(x,t)={-\gamma(x)v_t(x,t)}\in L^1([0,T],L^2)$ on the time interval
$[0,T]$,
we can go from $V(t_0)$
to $V(t_0+T)$ in the control problem \eqref{eq}.

\medskip

\noindent {\bf Basic fact 3:} if $(v,v_t)(t)$ is a solution of the damped
wave equation \eqref{eq-damped}, then, by using the control
$u(x,t)={\gamma(x)v_t(x,t_0-t)}$ on the time interval $[0,T]$, we can go
from
$(v,-v_t)(t_0)$ to $(v,-v_t)(t_0-T)$ in the control problem \eqref{eq}.

\medskip

\noindent {\bf Step 4: the ``double U-turn'' argument.} \\
This step consists in proving the following property.
\begin{prop}If
$V(t)=(v,v_t)(t)$ is a globally bounded solution of the damped wave equation
\eqref{eq-damped} (i.e. $V(t)$ belongs to the attractor $\Ac$), then for any
points $V(t_0)$ and $V(t_1)$ of the trajectory, there exists a control $u$ on a
time interval $[0,T]$ bringing $V(t_0)$ to $V(t_1)$.
\end{prop} 
\begin{demo}
If $t_1\geq t_0$, this is simply Basic Fact 2. Assume that $t_1<t_0$, we proceed
as follows. By Theorem \ref{th-dyn}, the damped wave equation \eqref{eq-damped} generates
a gradient dynamical system. Thus, there exist two equilibrium points $e_\pm$
and two times
$t_\pm$ with $t_-<t_1<t_0<t_+$ such that $V(t_\pm)$ belongs to $\Nc(e_\pm)$. We
start from $V(t_0)$ and using Basic Fact 2 with the control {$u=-\gamma
v_t(t_0+t)$}, we reach $V(t_+)=(v,v_t)(t_+)$ which belongs to $\Nc(e_+)$. Notice
that $\Nc(e_+)$ is assumed to be a ball centered in $(e_+,0)$, and thus that
$(v,-v_t)(t_+)$ also belongs to $\Nc(e_+)$. Due to Basic Fact 1, we can find a
control to go from $(v,v_t)(t_+)$ to $(v,-v_t)(t_+)$. Then, Basic Fact 3 shows
that we can travel backward to reach $(v,-v_t)(t_-)$. Using Basic Fact 1 again,
we go from $(v,-v_t)(t_-)$ to $V(t_-)=(v,v_t)(t_-)$. Finally, we simply follow
the trajectory as described in Basic Fact 2 to reach $V(t_1)$.
\end{demo}

\medskip

\noindent {\bf Step 5: the reachable set $\mathcal{R}(e,0)$ of an equilibrium $(e,0)$ is open in $X$.} Let
$e$ be an equilibrium point of the wave equation. Assume that there exists a
control $u$ on a time interval $[0,T]$ such that the solution of \eqref{eq}
satisfies $V(0)=(e,0)$ and $V(T)=V_1$. By Theorem \ref{th-Cauchy1} and the
reversibility of the equation, we know that if a control $u\in L^1([0,T],L^2)$
is fixed, the flow map $\Phi$ which sends an initial data $V_0$ to a final data
$V(T)$ at time $T$ for a solution of \eqref{eq} is an homeomorphism of $X$. In
particular, there exists a neighborhood $\Uc$ of 
$V_1$ such that the backward Cauchy problem \eqref{eq} starting at $t=T$ in
$\Uc$ with control $u$ is defined in $[0,T]$ and arrives in the
neighborhood $\Nc(e)$ of $(e,0)$. In other words, applying 
the control $u$ in \eqref{eq}, we can reach any point of $\Uc$ by starting
from some $V_0\in\Nc(e)$ and applying the control $u$. 
On the other hand, due to the local control hypothesis
(Basic Fact 1), there exist controls $\tilde u$ bringing $(e,0)$ to any point
of the ball $\Nc(e)$. Therefore, applying successively the controls $\tilde u$
and $u$, we can reach any point of $\Uc$ from $(e,0)$.       

\medskip

\noindent {\bf Step 6: the reachable set $\mathcal{R}(e,0)$ of an equilibrium $(e,0)$ is closed in the
attractor $\Ac$.} Let $V_1^n$ be points of the compact global
attractor $\Ac$ such that there exist controls $u^n$ on $[0,T^n]$ bringing
$(e,0)$ to $V_1^n$ and such that $(V_1^n)$ converges to $V_1\in\Ac$. We denote
by $V(t)$ the solution of the damped wave equation \eqref{eq-damped} with
$V(0)=V_1$. Since the flow $S(t)$ of \eqref{eq-damped} is assumed to be
gradient, there exists $T>0$ and $\widetilde{e}$ an equilibrium point such
that $V(T)$ belongs to $\Nc(\widetilde{e})$. By continuity of the flow of
\eqref{eq-damped},
$S(T)V_1^n$ also belongs to $\Nc(\widetilde{e})$ for $n$ large enough.
Using successively Basic Fact 2, Basic Fact 1 and Step 4, we can go from $V_1^n$
to $V(0)=V_1$ via $S(T)V_1^n$ and $V(T)$.

\medskip

\noindent {\bf Step 7: conclusion.} The compact global attractor $\Ac$ of
$S(t)$ is a connected set (see Theorem \ref{th-dyn}).
Therefore, Steps 5 and 6 show that the
reachable set of an equilibrium point is a neighborhood $\Nc(\Ac)$ of the
attractor $\Ac$. In particular, we can go from the neighborhood of any
equilibrium point to the one of any other equilibrium point. Let $V_0$
and $V_1$ be two points of $X$. Let $V(t)$ and $\tilde V(t)$ be the 
trajectories of \eqref{eq-damped} satisfying $V(0)=V_0$ and $(\tilde v,-\tilde
v_t)(0)=V_1$. For $T$ large enough, $V(T)$ and $\tilde V(T)$ are in
neighborhoods of equilibrium points, and so is $(\tilde v,-\tilde v_t)(T)$.
Combining the previous argument, we can go from $V_0$ to $V_1$ through the
control problem \eqref{eq}, via $V(T)$ and $(\tilde
v,-\tilde v_t)(T)$.

\medskip

\noindent {\bf Step 8: uniformity of the time of control.} By assumption, the
time of local control is uniformly bounded by $T(e)$ in a neighborhood $\Nc(e)$
of an equilibrium point. Since the set of all the equilibria is a closed subset
of the compact attractor $\Ac$ and is therefore compact, the time of control
$T(e)$ can be chosen to be independent of $e$ and we can consider only a
finite number of equilibrium points in the above arguments. For any
neighborhood $\Nc(\Ac)$ of $\Ac$ and any ball $B_X(0,R)$, there is a time $T$
such that $S(T)B_X(0,R)\subset \Nc(\Ac)$. Considering the paths yielded by the
above arguments to link two given points of $B_X(0,R)$, there is only one last
property to show: there is a time $T(\Ac)$ and a neighborhood $\Nc(\Ac)$ of $\Ac$ such that, for any $V_0\in \Nc(\Ac)$,
$S(t)V_0$ belongs to a neighborhood $\Nc(e)$ for some $t\in [0,T(\Ac)]$. To show
this last property, we argue by contradiction. Assume that there exist
sequences $(V^n_0)\subset X$, with $d(V^n_0,\Ac)\leq 1/n$, and $(T^n)\rightarrow +\infty$ such that
$S(t)V_0^n$ does not cross any neighborhood $\Nc(e)$ for $t\in [0,T^n]$. Since
$\Ac$ is compact, we can assume that $(V^n_0)$ converges to $V_0\in\Ac$.
Because $S(t)$ is a gradient dynamical system, $S(T)V_0$ belongs to a
neighborhood $\Nc(e)$ for $T$ large enough, and thus $S(T)V_0^n$ also belongs
to $\Nc(e)$ for $n$ large enough. This contradicts the definition of $T^n$ and
the fact that $T^n$ goes to $+\infty$.


\section{Discussion}\label{sect-disc}

{First notice that Theorem \ref{th} is a result of semi-global
controllability in the sense that the time of control depends on the sizes of
the initial and final data. One may expect that, for the damped wave equation
\eqref{eq-damped}, the large balls converge, with a uniform exponential rate,
to an absorbing ball $\Bc_0$, which contains the attractor $\Ac$. This would
implies that the time of control in a large ball $\Bc$ of radius $R$ is of order
$T(\Bc)\sim C\ln R+T(\Bc_0)$. However, proving such a uniform exponential rate
of convergence to the absorbing ball is a difficult problem related to the
uniform stabilization of the semilinear wave equation (see
\cite{Zuazua-stab} and its application to the convergence to the absorbing ball
in Proposition 3.6 of \cite{RJ}).}

\begin{figure}[p]
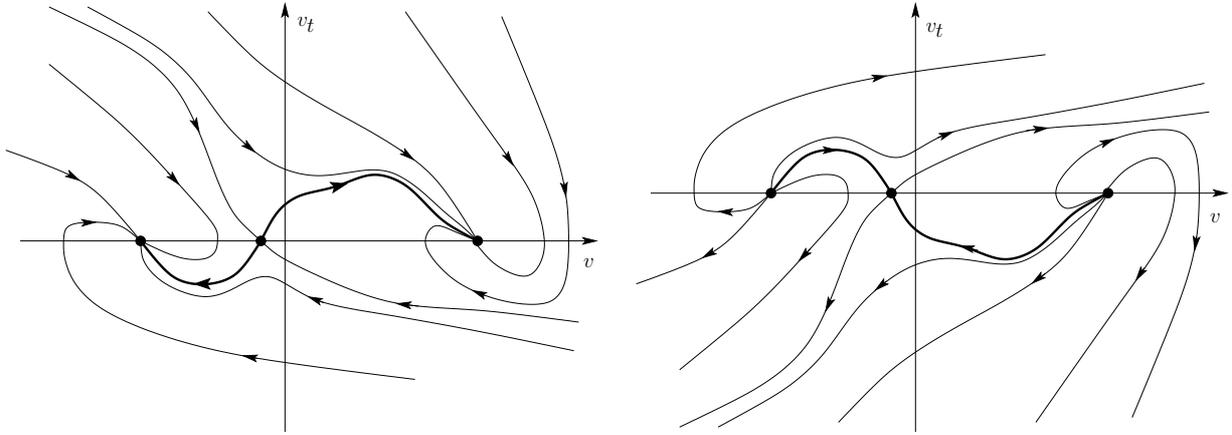

\begin{center}
\resizebox{1\textwidth}{!}{\input{flot-pos.pstex_t}~~~~~~~~~
\input{flot-neg.pstex_t } }
\end{center}
\caption{\it Right: the flow of the dynamical system $S(t)$ generated by the
feedback control {$u=-\gamma(x) v_t$}. The flow is represented in the phase
plane
$(v,v_t)$. The compact global attractor $\Ac$ (in bold) consists in three
equilibrium points and two heteroclinic orbits connecting them.
Left: the flow generated by the feedback control {$u=\gamma(x) v_t$}. It
is deduced from the flow of $S(t)$ by reversing time and orientation of the
second coordinate.}
\label{fig1}
\end{figure}

\begin{figure}[p]
\begin{center}
\resizebox{0.9\textwidth}{!}{\input{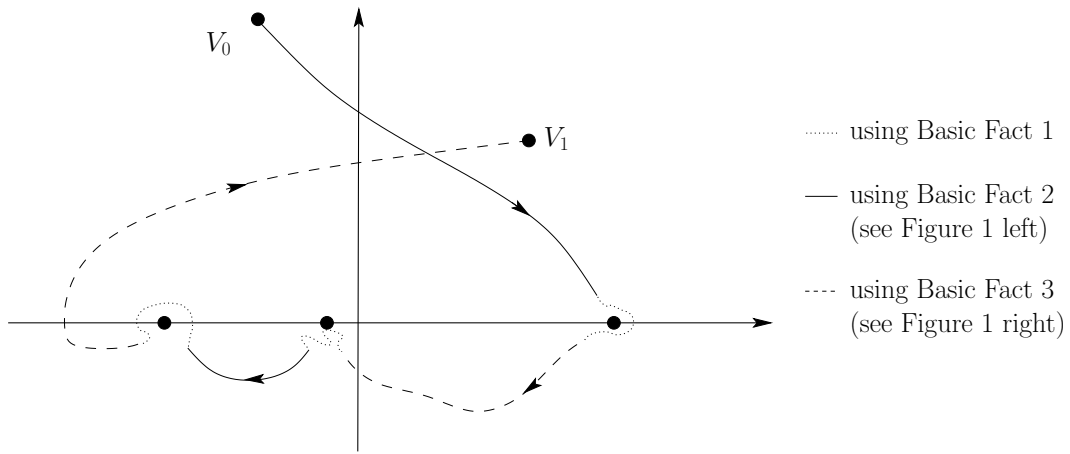}}
\end{center}
\caption{\it An example of global control using successively the Basic Facts
1, 2 and 3 introduced at the beginning of the proof of Theorem \ref{th}. The
resulting trajectory consists in switching from one of the flows of Figure
\ref{fig1} to the other, using in between the local controls near the
equilibrium points.}
\label{fig3}
\end{figure}

The arguments of the proof of Theorem \ref{th} are very related to the
qualitative dynamics of the wave equation \eqref{eq} with the feedback control
{$u=-\gamma(x) v_t$}. In the simplest cases, one can follow the arguments
to
construct a simple control as illustrated in Figures \ref{fig1} and
\ref{fig3}.

Notice that the proof of Theorem \ref{th} is somehow
constructive: the way to
travel between two equilibrium points is described in 
terms of the heteroclinic orbits and moreover the local control may be explicit
if it is obtained via a Banach fixed point theorem. In this sense, our method of proof is more
explicit than the Leray-Schauder fixed point argument used in the proofs of
several control results for the nonlinear wave equations (see
\cite{Zuazua-control} for example). Of course, our proof of Theorem \ref{th} is
not really constructive, because the computation of heteroclinic orbits is not
explicit. But it is reasonable to expect some cases where approximations of the
heteroclinic orbits are numerically known. In this case, our proof of Theorem
\ref{th} gives an explicit way to compute an approximated control to connect
two equilibrium points. By the way, the question whether Theorem \ref{th} can be
proved by using
a fixed point method as the one of \cite{Zuazua-control}, or not, is an open
interesting question. In any case, we think that our method of proof is new and 
represents an interesting alternative to the proof by Leray-Schauder theorem,
which is actually available only for weak nonlinearities.

{As noticed after Theorem \ref{th-dyn}, it is
reasonable to think the attractor $\Ac$ of
\eqref{eq-damped} as a finite number of isolated equilibrium points connected
by heteroclinic orbits, as in Figure \ref{fig1}.}
In this case, one can go from any equilibrium
point to another one by following a finite number of heteroclinic orbits,
forward or backward in time. In this
point of view, the steps 5 and 6 of the proof of Theorem \ref{th} may seem
unnatural. However, the above generic framework is
not true in general: $\Ac$ may contain a continuum of equilibrium points, which
precludes the connection between two parts of $\Ac$ via 
heteroclinic orbits. That is why, in the proof of Theorem \ref{th}, we used the
topological connectedness of $\Ac$ and not the connectedness in the sense of
path of heteroclinic orbits.

\vspace{10mm}

\noindent {\bf Acknowledgements:} the authors are grateful to Emmanuel Tr\'elat
and Enrique Zuazua for interesting discussions. 
{They also thank the referees for their remarks and corrections.}
The work of the second author has been supported by ANR EMAQS.



\end{document}